\documentclass[12pt]{article}
\def\boxit{$\sqcap\kern-8pt\sqcup$}
\def\littbox{\hfill{\boxit{}}\vskip .2cm}

\newtheorem{teo}{Theorem}[section]
\newtheorem{lem}[teo]{Lemma}
\newtheorem{cor}[teo]{Corollary}

\newtheorem{rem}[teo]{Remark}

\newtheorem{prob}[teo]{Problem}

\newcommand{\nene}{{\rm I}\!{\rm N}}

\input epsf

\begin{document}

\title{On the Frobenius Number of Fibonacci Numerical Semigroups}
\author{J.M. Mar\'{\i}n\thanks{{\em e-mail}: jmarin@us.es}\\
{\it Dep. Matem\'{a}tica Aplicada I}\\ 
{\it Universidad de Sevilla. Avda. Reina Mercedes s/n. 41012}\\
{\it Sevilla, Spain}\\
J. L. Ram\'{\i}rez Alfons\'{\i}n\thanks{{\em e-mail}: ramirez@math.jussieu}\\
{\it Equipe Combinatoire et Optimisation,} \\
{\it Universit\'{e} Pierre et Marie Curie (Paris 6),}\\
{\it Case 189 - 4 Place Jussieu, 75252 Paris, Cedex 05,  France}\\
and \\
M.P. Revuelta\thanks{{\em e-mail}: pastora@us.es}\\
{\it Dep. Matem\'{a}tica Aplicada I}\\ 
{\it Universidad de Sevilla. Avda. Reina Mercedes s/n. 41012}\\
{\it Sevilla, Spain}}

\maketitle

\begin{abstract} In this paper we compute the Frobenius number of  certain {\em Fibonacci numerical semigroups}, that is,
numerical semigroups generated by a set of Fibonacci numbers, in terms of Fibonacci numbers. 
\end{abstract}

\section{Introduction}

Let $s_1,\dots ,s_n$ be positive integers such that their greatest commun divisor is one.
Let $S=<s_1,\dots ,s_n>$ be the numerical semigroup\footnote{Recall that a {\em semigroup } $(S,*)$ consists of a nonempty set 
$S$ and an associative binary operation $*$ on $S$. If, in addition, there exists an element, which is usually 
denoted by $0$, in $S$ such that $a+0=0+a=a$ for all $a\in S$, we say that $(S,*)$ is a {\em monoid}.
A {\em numerical semigroup} is a submonoid of $\nene$ such that the greatest common divisor of 
its elements is equal to one.} generated by $s_1,\dots ,s_n$. A {\em Fibonacci numerical semigroup} is a numerical semigroup generated 
by a set of Fibonacci numbers $F_{i_1},\dots ,F_{i_r}$, for some integers $3\le {i_1}<\cdots <{i_r}$ where $g.c.d.(F_{i_1},\dots ,F_{i_r})=1$.
\vskip .3cm

The so-called {\em Frobenius number}, denoted by $g(s_1,\dots ,s_n)$, is  defined as the largest 
integer not belonging to $S$, that is, the largest integer that is not
representable as a nonnegative integer combination of $s_1,\dots ,s_n$.
It is well known that $g(s_1,s_2)=s_1s_2-s_1-s_2$.
In general, finding $g(S)$ is a difficult problem and so formulas and upper bounds for particular sequences are of interest.
For instance, it is known \cite{Ro1} the value of $g(S)$ when $S$ is an arithmetical sequence

\begin{equation}\label{arit}
g(a,a+d,\dots ,a+kd)=a\left(\left\lfloor{{a-2}\over {k}}\right\rfloor\right)+d(a-1)
\end{equation}

We refer the reader to \cite{Ram4} where a complete account on the Frobenius problem can be found.
\vskip .3cm

In this paper, we investigate the value of $g(F_i,F_j,F_l)$ for some triples $3\le i<j<l$ (we always assume that $g.c.d.(F_i,F_j,F_l)=1$, 
recall that $g.c.d.(F_i,F_{i+l})=1$ if $i\not\hskip-0.10em{\mid}\hskip0.30em l$).
\vskip .3cm

We first notice that $g(F_i,F_{i+1},F_{l})=g(F_i,F_{i+1})$ for any integer $l\ge i+2$. Indeed, since 
$F_l=F_{i+m}=F_mF_{i+1}+F_{m-1}F_i$ is a nonnegative integer combination of $F_i$ and $F_{i+1}$ then the semigroups 
$<F_i,F_{i+1},F_l>$ and $<F_i,F_{i+1}>$ generate the same set of elements and thus they have the same Frobenius number.  
\vskip .3cm

Let us consider then $g(F_i,F_{i+2},F_{l})$ with $l\ge i+3$. We notice that the case when $l=i+3$
is a consequence of equation (\ref{arit}) since the triple $\{F_i,F_{i+2},F_{i+3}\}=\{F_i,F_i+F_{i+1},F_i+2F_{i+1}\}$ form an arithmetical sequence.
However, it can be checked that $\{F_i,F_{i+2},F_{i+k}\}$ do not form an arithmetical sequence when $k\ge 3$ and the calculation of
$g(F_i,F_{i+2},F_{i+k})$ is more complicated. 

We state our main result.

\begin{teo}\label{theo2} Let $i,k\ge 3$ be integers and let $r=\lfloor{{F_i-1}\over {F_k}}\rfloor$. Then, 

$$g(F_{i},F_{i+2},F_{i+k})=\left\{\begin{array}{ll}
(F_i-1)F_{i+2}-F_i(rF_{k-2}+1)& \hbox{if $r=0$ or $r\ge 1$ and}\\
& \hbox{$F_{k-2}F_i<(F_i-rF_k)F_{i+2}$,}\\
\\
(rF_k-1)F_{i+2}-F_i((r-1)F_{k-2}+1)& \hbox{otherwise.}\\
\end{array}\right.$$

\end{teo}

Let $N(a_1,\dots ,a_n)$ be the number of positive integers with
no representation by a nonnegative integer combination of $a_1,\dots ,a_n$. 
We refer the reader to \cite[Chapter 5]{Ram4} for results related to $N(a_1,\dots ,a_n)$.
Theorem \ref{theo2} yields to the following result.

\begin{cor}\label{corr} Let $i,k\ge 3$ be integers and let $r=\lfloor{{F_i-1}\over {F_k}}\rfloor$.Then, 
$$N(F_{i},F_{i+2},F_{i+k})={{(F_i-1)(F_{i+2}-1)-rF_{k-2}(2F_i-F_k(1+r))}\over {2}}\cdot$$
\end{cor}

\section{Fibonacci semigroups}

In order to prove Theorem \ref{theo2} we need the following result due to Brauer and Shockley \cite{BS}.

\begin{lem}\label{l1} Let $1<a_1<\cdots <a_n$ be integers with $g.c.d.(a_1,\dots ,a_n)=1$. Then,
$$g(a_1,\dots ,a_n)=
\max\limits_{l\in \{1,2,\dots ,a_n-1\}}\{t_l\}-a_1$$
where $t_l$ is the smallest positive integer congruent to $l$ modulo $a_1$,
that is representable as a nonnegative integer combination of
$a_2,\dots ,a_{n}$.
\end{lem}

{\em Proof.} Let $L$ be a positive integer. If
$L\equiv 0 \bmod {a_1}$ then $L$ is a nonnegative integer combination of $a_1$.
If $L\equiv l \bmod {a_1}$ then $L$ is a nonnegative integer combination of
$a_1,\dots ,a_n$ if and only if $L\ge t_l$. \littbox
\vskip .3cm

{\em Proof of Theorem \ref{theo2}.} Let $T^*=\{t^*_0,\dots ,t^*_{F_i-1}\}$ where $t^*_l$ be the smallest positive integer congruent to $l$ modulo $F_i$,
that is representable as a nonnegative integer combination of
$F_{i+2}$ and $F_{i+k}$. We shall find $t^*_l$ for each $l=0,1,\dots ,F_i-1$.
To this end, we consider all nonnegative integer combinations of $F_{i+2}$ and $F_{i+k}$. We construct the following table, denoted by $T_1$, having as entry $t_{x,y}$ the combination 
of the form $xF_{i+2}+yF_{i+k}$ with integers $x,y\ge 0$, see below.

$$\begin{array}{c|ccccc}
x\backslash y & 0 & 1 & 2 & &\cdots\\
\hline
0 & 0 & F_{i+k} & 2F_{i+k} & &\cdots \\
1 & F_{i+2} &   F_{i+k}+F_{i+2} & 2F_{i+k}+F_{i+2} &&\cdots \\
2 &  2F_{i+2} & F_{i+k}+2F_{i+2} & 2F_{i+k}+2F_{i+2} &&\cdots \\
3 &3F_{i+2} &  F_{i+k}+3F_{i+2} &2F_{i+k}+3F_{i+2}  &&\cdots\\
\vdots & \vdots & \vdots & \vdots& \\
F_k-1 & (F_k-1)F_{i+2}& F_{i+k}+(F_k-1)F_{i+2}&2F_{i+k}+(F_k-1)F_{i+2}&&\cdots\\
\vdots & \vdots & \vdots & \vdots\\
\end{array}$$

We notice that 
$$F_{i+k}=F_{k-2}F_{i+1}+F_{k-1}F_{i+2}=F_{k-2}(F_{i+2}-F_i)+F_{k-1}F_{i+2}=F_{i+2}F_k-F_{k-2}F_i$$ 

so, we obtain that

$$xF_{i+2}+yF_{i+k}=xF_{i+2}+y(F_{i+2}F_k-F_{k-2}F_i)=(x+yF_k)F_{i+2}-yF_{k-2}F_i.$$
 
Thus, $T_1$ can also be given by the following table, denoted by $T_2$,
 
{\scriptsize
$$\label{tableb}\begin{array}{c|ccccccc}
x\backslash y & 0 & 1 & 2 & \cdots & r & \cdots\\
\hline
0 & 0 & F_kF_{i+2}-F_{k-2}F_i  & 2F_kF_{i+2}-2F_{k-2}F_i& \cdots & rF_kF_{i+2}-rF_{k-2}F_i& \cdots\\
1 & F_{i+2} & (1+F_k)F_{i+2}-F_{k-2}F_i   & (1+2F_k)F_{i+2}-2F_{k-2}F_i &\cdots & (1+rF_k)F_{i+2}-rF_{k-2}F_i&\cdots \\
2 &  2F_{i+2} & (2+F_k)F_{i+2}-F_{k-2}F_i & (2+2F_k)F_{i+2}-2F_{k-2}F_i &\cdots & (2+rF_k)F_{i+2}-rF_{k-2}F_i&\cdots\\
\vdots & \vdots & \vdots & \vdots&  &\vdots\\
l &lF_{i+2} &  (l+F_k)F_{i+2}-F_{k-2}F_i& (l+2F_k)F_{i+2}-2F_{k-2}F_i&\cdots &(l+rF_k)F_{i+2}-rF_{k-2}F_i&\cdots\\
\vdots & \vdots & \vdots & \vdots& & \vdots\\
F_k-1 & (F_k-1)F_{i+2}&(2F_k-1)F_{i+2}-F_{k-2}F_i &(3F_k-1)F_{i+2}-2F_{k-2}F_i&&\cdots\\
\vdots & \vdots & \vdots & \vdots&& \vdots\\
\end{array}$$}

Let $S$ be the set formed by the first $F_k-1$ entries of columns zero, one, two, and so on, that is,

$$S=\{t_{0,0},t_{1,0},\dots ,t_{F_{k}-1,0},t_{0,1},t_{1,1},\dots ,t_{F_{k}-1,1},\dots, t_{0,r},t_{1,r},\dots ,t_{F_{k}-1,r},\dots\}.$$

\begin{rem}\label{rree}

\begin{enumerate}

\item [(a)] Let $r=\lfloor{{F_i-1}\over {F_k}}\rfloor$ and set $F_i-1=rF_k+l$ for some integer $0\le l\le F_k-1$. Let 

$$S'=\{t_{0,0},t_{1,0},\dots ,t_{F_{k}-1,0},t_{0,1},t_{1,1},\dots ,t_{F_{k}-1,1},\dots, t_{2,r},t_{1,r},\dots ,t_{l,r}\},$$

Then, for each $t_{x,y}=(x+yF_k)F_{i+2}-yF_{k-2}F_i \in S'$ we have that $0\le x+yF_k\le F_i-1$. Moreover, since $g.c.d.(F_{i +2},F_{i})=1$ then 
$S'$ forms a complete system of rests modulo $F_i$. 

\item [(b)] The elements of $S$ can be represented as $s_{x}=xF_{i+2}-\lfloor {x\over  {F_k}} \rfloor F_{k-2}F_{i}$ for $x=0,1,\dots$. Indeed, it can be checked
that $S=\bigcup_{q\ge 1}S_q$ where

$$S_q=\{s_{qF_k},s_{qF_k+1},\dots ,s_{(q+1)F_k-1}\}=\{t_{0,q},\dots ,t_{F_k-1,q}\}$$ 

for each integer $q=0,1,2,\dots$.

\item [(c)] By using table $T_2$ we have that $t_{i,j}<t_{k,l}$ for all $i\le k$ and all $j\le l$.
\end{enumerate}
\end{rem}

\begin{lem}
Let $t_{u,v}$ be an entry of $T_1$ such that $t_{u,v}\not\in S'$. There exists $t_{x,y}\in S'$ such that $t_{u,v}\equiv t_{x,y} \bmod F_i$ and $t_{u,v}>t_{x,y}$.
\end{lem}

{\em Proof.} We first notice that the set $S$ can be written as follows

$$\begin{array}{ccccccccccc}
\{s_0,& \dots &,s_{F_k-1},&s_{F_k},&\dots &,s_{2F_k-1},&\dots &,s_{rF_k},&\dots &,s_{rF_k+l}=s_{F_i-1},\\
s_{F_i},& \dots &,s_{F_i+F_k-1},&s_{F_i+F_k},&\dots &,s_{F_i+2F_k-1},&\dots &,s_{F_i+rF_k},&\dots &,s_{2F_i-1},\\
s_{2F_i},& \dots &,s_{2F_i+F_k-1},&s_{2F_i+F_k},&\dots &,s_{2F_i+2F_k-1},&\dots &,s_{2F_i+rF_k},&\dots &,s_{3F_i-1},\dots\}\\
\end{array}$$

where $S'=\{s_0,\dots ,s_{F_k-1},s_{F_k},\dots ,s_{2F_k-1},\dots ,s_{rF_k},\dots ,s_{F_i-1}\}$. We have two cases.
\vskip .3cm

Case A) Suppose that $t_{u,v}\in S\setminus S'$. Then $t_{u,v}$ is of the form $s_{pF_i+g}$ for some integers $p\ge 1$ and $0\le g\le F_i-1$.
It is clear that, 

$$s_g=gF_{i+2}-\left\lfloor{g\over {F_k}}\right\rfloor F_iF_{k-2}\equiv (pF_i+g)F_{i+2}-\left\lfloor{{pF_i+g}\over {F_k}}\right\rfloor F_iF_{k-2}=g_{pF_i+g}\bmod F_i.$$

We will show that $s_{pF_i+g}>s_g$. To this end, it suffice to prove that $s_{F_i+g}>s_g$ (since $s_{pF_i+g}\ge s_{F_i+g}$). 
Recall that $r=\lfloor{{F_i-1}\over {F_k}}\rfloor$ and that $F_i-1=rF_k+l$ for some integer $0\le l\le F_k-1$. We have two subcases.
\vskip .3cm

Subcase a) If $r=0$ then $F_k\ge F_i$.  If $F_k=F_i$ then $s_{F_i+g}=t_{g,1}$ and, by Remark \ref{rree} (c), $t_{g,0}<t_{g,1}$.
If $F_k>F_i$ then $s_{F_i+g}=t_{q,0}$ for some integer $q\ge F_i$ and, by Remark \ref{rree} (c), $t_{g,0}<t_{q,0}$.
\vskip .3cm

Subcase b) If $r\ge 1$ then $s_{F_i+g}>s_g$ holds if and only if 

$$(F_i+g)F_{i+2}-\left\lfloor{{F_i+g}\over {F_k}}\right\rfloor F_iF_{k-2}>gF_{i+2}-\left\lfloor{g\over {F_k}}\right\rfloor F_iF_{k-2}$$
or equivalently if and only if
$$F_{i+2}>F_{k-2}\left( \left\lfloor{{F_i+g}\over {F_k}}\right\rfloor - \left\lfloor{g\over {F_k}}\right\rfloor\right)\cdot$$

Let $g=mF_k+n$ with $0\le n\le F_k-1$. Since $F_i-1=rF_k+l$ with $0\le l\le F_k-1$ then 

$$\left\lfloor{{F_i-1+g+1}\over {F_k}}\right\rfloor=\left\lfloor{{rF_k+l+mF_k+n+1}\over {F_k}}\right\rfloor\le r+m+1$$

and thus

$$\left\lfloor{{F_i+g}\over {F_k}}\right\rfloor - \left\lfloor{g\over {F_k}}\right\rfloor\le r+m+1-m=r+1.$$

So, it is enough to show that $F_{i+2}>(r+1)F_{k-2}$ or equivalently to show that $F_i+F_{i+1}>(r+1)F_{k-2}$. Since $F_i=rF_k+l+1$ then the latter inequality holds if and only if
$rF_k+l+1+F_{i+1}>rF_{k-2}+F_{k-2}$, that is, if and only if 
$$r(F_k-F_{k-2})+l+1+F_{i+1}=r(F_{k-1})+l+1+F_{i+1}>F_{k-2}$$ 
which is true since $r\ge 1$.
\vskip .3cm

Case B)  Suppose that $t_{u,v}\not\in S$. Then we have that $0\le x\le F_k-1<u$.  If $v\ge y$ then, by Remark \ref{rree} (c), $t_{x,y}<t_{x,v}<t_{u,v}$.
So, we suppose that $v<y$. Since, $t_{u,v}\equiv t_{x,y} \bmod F_i$ then $u+vF_k\equiv x+yF_k \bmod F_i$ but, by Remark \ref{rree} (a), $0\le x+yF_k\le F_i-1$
so $u+vF_k=d(x+yF_k)$ for some integer $d\ge 1$ and thus $u+vF_k\ge x+yF_k$. Also, since $v<y$, then $-vF_{k-2}F_i>-yF_{k-2}F_i$. So, 
combining the last two inequalities we have that $$t_{u,v}=(u+vF_k)F_{i+2}-vF_{k-2}F_i> (x+yF_k)F_{i+2}-yF_{k-2}F_i=t_{x,y}.$$
\littbox

Therefore, by the above lemma, we have that for each $x=0,\dots ,F_i-1$, $s_x$ is the smallest 
positive integer congruent to $l$ modulo $F_i$, for some integer $0\le l\le F_i-1$,
that is representable as a nonnegative integer combination of $F_{i+2}$ and $F_{i+k}$, that is, $S'=T^*$. 

Now, by Remark \ref{rree} (c), if $r\ge 1$ then

$$\hbox{ $\ \  t_{F_k-1,i}=\max\limits_{0\le x\le F_k-1}\{t_{x,i} | t_{x,i}\in S'\}$ for each $i=0,\dots ,r-1$,}$$

$$t_{F_k-1,r-1}=\max\limits_{0\le i\le r-1}\{t_{F_k-1,i} | t_{F_k-1,i}\in S'\},$$

and

$$t_{l,r}=\max\limits_{0\le x\le l}\{t_{x,r} | t_{x,r}\in S'\}.$$

Thus, $$\max\{s | s\in S'\}=\left\{\begin{array}{ll}
t_{l,r}& \hbox{if $r=0$,}\\
\max\{t_{F_k-1,r-1},t_{l,r}\}&\hbox{otherwise.}
\end{array}\right.$$

The result follows since $t_{l,r}>t_{F_k-1,r-1}$ if and only if
$$(rF_k+l)F_{i+2}-rF_{k-2}F_i=(F_i-1)F_{i+2}-rF_{k-2}F_i>(rF_k-1)F_{i+2}-(r-1)F_{k-2}F_i$$

or equivalently, if and only if $F_{i+2}(F_i-rF_k)>F_{k-2}F_i$. \littbox

We will use the following result due to Selmer \cite{Sel1} to show Corollary \ref{corr}.

\begin{lem}\label{sel} Let $1<a_1<\cdots <a_n$ be integers with $g.c.d.(a_1,\dots ,a_n)=1$.
If $L=\{1,\dots ,a_1-1\}$ then

$$N(a_1,\dots ,a_n)={1\over {a_1}}\sum\limits_{l\in L}t_l-{{a_1-1}\over 2}\cdot$$

where $t_l$ is the smallest positive integer congruent to $l$ modulo $a_1$,
that is representable as a nonnegative integer combination of
$a_2,\dots ,a_{n}$.
\end{lem}

{\em Proof.} The number of $M\equiv l\not\equiv 0\bmod {a_1}$ with $0<M<t_l$
is given by $\lfloor {{t_1}\over {a_1}} \rfloor$. By assuming that $0<l<a_1$, we
have $\lfloor {{t_l}\over {a_1}} \rfloor= {{t_l-l}\over {a_1}}$. The result
follows by summing over $l\in L$. \littbox
\vskip .3cm

{\em Proof of Corollary \ref{corr}.} Let $r=\lfloor{{F_i-1}\over {F_k}}\rfloor$ and set $F_i-1=rF_k+l$ for some integer $0\le l\le F_k-1$.
By Lemma \ref{sel} and Remark \ref{rree} (b), we have

$$\begin{array}{ll}
N(F_i,F_{i+2},F_{i+k})& ={1\over {F_i}}\sum\limits_{s\in S'}s-{{F_i-1}\over 2}\\
&={1\over {F_i}}\sum\limits_{j=0}^{F_i-1}(jF_{i+2}-F_{k-2}\lfloor {j\over  {F_k}} \rfloor F_{i})-{{F_i-1}\over 2}\\
&={1\over {F_i}}\left(F_{i+2}{{(F_i-1)F_i}\over 2}\right)-{1\over {F_i}}(F_{k-2}F_i)\sum\limits_{j=0}^{F_i-1}\left\lfloor {j\over  {F_k}} \right\rfloor -{{F_i-1}\over 2}.\\
\end{array}$$

By using the table $T_1$, it is easy to verify that 

$$\sum\limits_{j=0}^{F_i-1}\left\lfloor {j\over  {F_k}} \right\rfloor= 0+F_k+2F_k+\cdots +(r-1)F_k+r(l+1)={{F_k(r-1)r}\over 2}+r(l+1)$$

and, since $l+1=F_i-rF_k$, then

$$\begin{array}{ll}
N(F_i,F_{i+2},F_{i+k})&={{F_{i+2}(F_i-1)}\over {2}}- F_{k-2}\left({{F_k(r-1)r}\over 2}+r(F_i-rF_k) \right)  -{{F_i-1}\over 2}\\
&={{(F_i-1)(F_{i+2}-1)}\over {2}}-F_{k-2}\left({{F_kr^2-F_kr+2F_ir-2r^2F_k}\over 2}\right)\\
&={{(F_i-1)(F_{i+2}-1)-rF_{k-2}(2F_i-F_k(1+r))}\over {2}}\cdot\\

\end{array}$$
\littbox
\vskip .3cm

We end with the following problem.

\begin{prob} Find upper (and lower) bounds (or formulas) for $g(F_i,F_j,F_k)$ for further triples $3\le i<j<k$. 
\end{prob}

{ AMS Classification number: 11B39, 11D04, 11A07}

\end{document}